\documentclass[fleqn,11pt,twoside]{article}

\usepackage{amsthm,amsthm,amssymb, color, xcolor,epsfig, graphics, subfigure}

\usepackage{amsmath, graphicx, latexsym, lscape }

\makeatletter
\newcommand{\copyrightnote}[2]{{\renewcommand{\thefootnote}{}
 \footnotetext{\small\it
\begin{flushleft}
 \copyright \ #1   #2
\end{flushleft}}}}

\newcommand{\Name}[1]{\begin{flushleft}
                       \LARGE \bf #1
                       \end{flushleft}\vspace{-3mm}}

\newcommand{\Author}[1]{\begin{flushleft}
                       \it #1 \end{flushleft}}

\newcommand{\Date}[1]{\begin{flushleft}
                      \small  \it #1 \end{flushleft}}

\newcommand{\evenhead}{Author \ name}
\newcommand{\oddhead}{Article \ name}

\renewcommand{\@evenhead}{
\hspace*{-3pt}\raisebox{-15pt}[\headheight][0pt]{\vbox{\hbox to \textwidth
{\thepage \hfil \evenhead}\vskip4pt \hrule}}}
\renewcommand{\@oddhead}{
\hspace*{-3pt}\raisebox{-15pt}[\headheight][0pt]{\vbox{\hbox to \textwidth
{\oddhead \hfil \thepage}\vskip4pt\hrule}}}
\renewcommand{\@evenfoot}{}
\renewcommand{\@oddfoot}{}

\long\def\@makecaption#1#2{%
  \vskip\abovecaptionskip
  \sbox\@tempboxa{\small \textbf{#1.}\ \ #2}%
  \ifdim \wd\@tempboxa >\hsize
    {\small \textbf{#1.}\ \ #2}\par
  \else
    \global \@minipagefalse
    \hb@xt@\hsize{\hfil\box\@tempboxa\hfil}%
  \fi
  \vskip\belowcaptionskip}

\newcommand{\JNMPnumberwithin}[3][\arabic]{%
  \@ifundefined{c@#2}{\@nocounterr{#2}}{%
    \@ifundefined{c@#3}{\@nocnterr{#3}}{%
      \@addtoreset{#2}{#3}%
      \@xp\xdef\csname the#2\endcsname{%
        \@xp\@nx\csname the#3\endcsname .\@nx#1{#2}}}}%
}

\newcommand{\resetfootnoterule} {
  \renewcommand\footnoterule{%
  \kern-3\p@
  \hrule\@width.4\columnwidth
  \kern2.6\p@}
}

\renewcommand{\footnoterule}{}

\makeatother

\theoremstyle{definition}

\begin{document}

\renewcommand{\evenhead}{ {\LARGE\textcolor{blue!10!black!40!green}{{\sf \ \ \ ]ocnmp[}}}\strut\hfill F. Calogero and F. Payandeh}
\renewcommand{\oddhead}{ {\LARGE\textcolor{blue!10!black!40!green}{{\sf ]ocnmp[}}}\ \ \ \ \   New solvable systems of 2 first-order nonlinearly-coupled ODEs}

\thispagestyle{empty}
\newcommand{\FistPageHead}[3]{
\begin{flushleft}
\raisebox{8mm}[0pt][0pt]
{\footnotesize \sf
\parbox{150mm}{{Open Communications in Nonlinear Mathematical Physics}\ \ \ {\LARGE\textcolor{blue!10!black!40!green}{]ocnmp[}}
\ \ Vol.2 (2022) pp
#2\hfill {\sc #3}}}\vspace{-13mm}
\end{flushleft}}

\FistPageHead{1}{\pageref{firstpage}--\pageref{lastpage}}{ \ \ Letter}

\strut\hfill

\strut\hfill

\copyrightnote{The author(s). Distributed under a Creative Commons Attribution 4.0 International License}

\qquad\qquad\qquad\qquad\qquad\qquad {\LARGE  {\sf Letter to the Editors}}

\strut\hfill

\Name{New solvable systems of 2 first-order nonlinearly-coupled ordinary differential equations}

\Author{Francesco Calogero$^{\,a,b,1,2}$ and Farrin Payandeh$^{\,c,3,4}$}

$^{a}$Physics Department, University of Rome "La Sapienza", Rome Italy

$^{b}$Istituto Nazionale di Fisica Nucleare, Sezione di Roma 1

$^{c}$Department of Physics, Payame Noor University (PNU), PO\ BOX

19395-3697 Tehran, Iran

$^{1}$francesco.calogero@uniroma1.it, $^{2}$francesco calogero@roma1.infn.it

$^{3}$f\_payandeh@pnu.ac.ir, $^{4}$farrinpayandeh@yahoo.com

\Date{Received September 14, 2022; Accepted September 26, 2022}

\setcounter{equation}{0}

\begin{abstract}
\noindent
In this short communication we introduce a rather simple
\textit{autonomous} system of $2$ \textit{nonlinearly-coupled} first-order
Ordinary Differential Equations (ODEs), whose \textit{initial-values}
problem is \textit{explicitly solvable} by \textit{algebraic} operations.
Its ODEs feature $2$ right-hand sides which are the ratios of $2$ \textit{%
homogeneous} polynomials of \textit{first} degree divided by the same
\textit{homogeneous} polynomial of \textit{second }degree. The model
features only $4$ \textit{arbitrary} parameters. We also report its \textit{%
isochronous} variant featuring $4$ \textit{nonlinearly-coupled} first-order
ODEs in $4$ dependent variables, featuring $9$ \textit{arbitrary}
parameters.\bigskip $\blacksquare $
\end{abstract}

\label{firstpage}



In this short communication we introduce a rather simple \textit{autonomous}
system of $2$ \textit{nonlinearly-coupled} first-order Ordinary Differential
Equations (ODEs) which features solutions obtainable by \textit{algebraic}
operations; and we report the explicit solution of its \textit{initial-values%
} problem. To the best of our knowledge these findings are new; but since
the corresponding literature includes an enormous number of entries over at
least $2$ centuries we cannot be quite certain.

The technique we have used to arrive at these results is fairly simple, but
we shall not describe it (the expert reader shall easily guess it from the
results reported below); we shall describe it in a separate paper \cite%
{CP2022} also reporting several other examples of \textit{algebraically
solvable} systems of $2$ first-order nonlinear ODEs.

The system of ODEs on which we focus reads as follows:
\begin{subequations}
\label{ODEs}
\begin{equation}
\dot{x}_{1}\left( t\right) =\frac{x_{1}\left( t\right) +\alpha
_{1}x_{2}\left( t\right) }{\beta _{1}\left[ x_{1}\left( t\right) \right]
^{2}+\left( \alpha _{1}\beta _{1}+\alpha _{2}\beta _{2}\right) x_{1}\left(
t\right) x_{2}\left( t\right) +\beta _{2}\left[ x_{2}\left( t\right) \right]
^{2}}~,
\end{equation}%
\begin{equation}
\dot{x}_{2}\left( t\right) =-\frac{x_{2}\left( t\right) +\alpha
_{2}x_{1}\left( t\right) }{\beta _{1}\left[ x_{1}\left( t\right) \right]
^{2}+\left( \alpha _{1}\beta _{1}+\alpha _{2}\beta _{2}\right) x_{1}\left(
t\right) x_{2}\left( t\right) +\beta _{2}\left[ x_{2}\left( t\right) \right]
^{2}}~.
\end{equation}

\textbf{Notation}: $t$ is the independent variable (one might think of it as
\textit{time}); $x_{n}\left( t\right) $ are the $2$ \textit{dependent}
variables; the superimposed dot denotes $t$-differentiation; $\alpha _{n}$
and $\beta _{n}$ ($n=1,2$) are $4$ ($t$-independent) \textit{arbitrary}
parameters. $2$ additional parameters $\gamma _{1}$ and $\gamma _{2}$ might
of course be introduced by rescaling the $2$ dependent variables ($%
x_{n}\left( t\right) \Rightarrow \gamma _{n}x_{n}\left( t\right) $); while
no third parameter $\gamma $ can be additionally introduced by rescaling the
independent variable $t$, since the system (\ref{ODEs}) is clearly \textit{%
invariant} under the common rescaling $x_{n}\left( t\right) \Rightarrow
\gamma x_{n}\left( t\right) $, $t\Rightarrow \gamma t$. $\blacksquare $

\textbf{Proposition}: the solution of the \textit{initial-values} problem of
this system (\ref{ODEs}) reads as follows:
\end{subequations}
\begin{equation}
x_{n}\left( t\right) =\left( \gamma _{n1}\sqrt{1+t/t_{1}}+\gamma _{n2}\sqrt{%
1+t/t_{2}}\right) ,~~~n=1,2~,  \label{x1x2}
\end{equation}%
with
\begin{subequations}
\label{aabbbb}
\begin{equation}
\gamma _{11}=\frac{b_{2}\left[ b_{1}x_{1}\left( 0\right) +a_{2}x_{2}\left(
0\right) \right] }{b_{1}b_{2}-a_{1}a_{2}}~,~~~\gamma _{22}=\frac{b_{1}\left[
a_{1}x_{1}\left( 0\right) +b_{2}x_{2}\left( 0\right) \right] }{%
b_{1}b_{2}-a_{1}a_{2}}~,
\end{equation}%
\begin{equation}
\gamma _{12}=\frac{-a_{2}\left[ a_{1}x_{1}\left( 0\right) +b_{2}x_{2}\left(
0\right) \right] }{b_{1}b_{2}-a_{1}a_{2}}~,~~~\gamma _{21}=\frac{-a_{1}\left[
b_{1}x_{1}\left( 0\right) +a_{2}x_{2}\left( 0\right) \right] }{%
b_{1}b_{2}-a_{1}a_{2}}~,
\end{equation}%
\begin{equation}
b_{1}=2\beta _{1}-\alpha _{2}\left( r+\alpha _{1}\beta _{1}+\alpha _{2}\beta
_{2}\right) ~,~~~b_{2}=-2\beta _{2}+\alpha _{1}\left( r+\alpha _{1}\beta
_{1}+\alpha _{2}\beta _{2}\right) ~,
\end{equation}%
~%
\begin{equation}
a_{n}=\left( -\right) ^{n+1}r+\alpha _{1}\beta _{1}-\alpha _{2}\beta
_{2}~,~~\ n=1,2~,
\end{equation}%
\begin{equation}
r=\sqrt{\left( \alpha _{1}\beta _{1}+\alpha _{2}\beta _{2}\right)
^{2}-4\beta _{1}\beta _{2}}~,  \label{r}
\end{equation}%
\begin{equation}
t_{1}=-\eta /\eta _{1}~,~~~t_{2}=\eta /\eta _{2}~,
\end{equation}%
\begin{equation}
\eta =\left( \gamma _{12}\gamma _{21}-\gamma _{11}\gamma _{22}\right)
\left\{ \beta _{1}\left[ x_{1}\left( 0\right) \right] ^{2}+\left( \alpha
_{1}\beta _{1}+\alpha _{2}\beta _{2}\right) x_{1}\left( 0\right) x_{2}\left(
0\right) +\beta _{2}\left[ x_{2}\left( 0\right) \right] ^{2}\right\} ~,
\end{equation}%
\begin{equation}
\eta _{1}=2\left[ \left( \alpha _{2}\gamma _{12}+\gamma _{22}\right)
x_{1}\left( 0\right) +\left( \gamma _{12}+\alpha _{1}\gamma _{22}\right)
x_{2}\left( 0\right) \right] ~,
\end{equation}%
\begin{equation}
\eta _{2}=2\left[ \left( \alpha _{2}\gamma _{11}+\gamma _{21}\right)
x_{1}\left( 0\right) +\left( \gamma _{11}+\alpha _{1}\gamma _{21}\right)
x_{2}\left( 0\right) \right] ~.~~~\blacksquare
\end{equation}

\textbf{Remark 1}. The formulas written above provide the \textit{explicit}
definition (up to the implicit sign ambiguities due to the square-roots
appearing in them) of the initial-values problem of the system of ODEs (\ref%
{ODEs}), in terms of the $4$ parameters $\alpha _{n}$ and $\beta _{n}$ ($%
n=1,2$) featured by the system of $2$ ODEs (\ref{ODEs}). The skeptical
reader who wishes to check that the formulas (\ref{x1x2}) with (\ref{aabbbb}%
) provide the solution of the initial-value problem of the system of ODEs (%
\ref{ODEs}) is of course welcome to do so! $\blacksquare $

\textbf{Remark 2}. Because of the sign ambiguities due to the square-roots
appearing in the formulas (\ref{x1x2}) and (\ref{r}), it might appear that
the formulas (\ref{x1x2}) are \textit{inadequate} to provide the solution of
the initial-values problem of the system of ODEs (\ref{ODEs}). But it is
easy to check that the formulas (\ref{x1x2}) with (\ref{aabbbb}) become
\textit{identities} at $t=0$, with the most obvious assignment of the
square-roots; and thereafter they provide a well-defined definition of the
solution for all time by \textit{continuity} in the independent variable $t$%
. Of course a \textit{singularity} may then be hit if the time-evolution
causes the argument of one of the $2$ square-roots $\sqrt{1+t/t_{n}}$ ($%
n=1,2 $) appearing in the right-hand of the $2$ eqs. (\ref{x1x2}) to vanish;
but this is a \textit{natural} feature of \textit{nonlinear} systems of
evolution equations. $\blacksquare $

\textbf{Remark 3}. Because the $2$ ODEs of the system (\ref{ODEs}) feature
right-hand sides which are \textit{homogeneous} in the $2$ dependent
variables $\tilde{x}_{1}\left( t\right) $ and $\tilde{x}_{2}\left( t\right) ,
$ and its solutions depend in a simple manner on the time variable $t$ (see (%
\ref{x1x2})), the more general system
\end{subequations}
\begin{subequations}
\label{isoiso}
\begin{equation}
\overset{\cdot }{\tilde{x}}_{1}\left( t\right) =\mathbf{i}\omega \tilde{x}%
_{1}\left( t\right) +\frac{\tilde{x}_{1}\left( t\right) +\alpha _{1}\tilde{x}%
_{2}\left( t\right) }{\beta _{1}\left[ \tilde{x}_{1}\left( t\right) \right]
^{2}+\left( \alpha _{1}\beta _{1}+\alpha _{2}\beta _{2}\right) \tilde{x}%
_{1}\left( t\right) \tilde{x}_{2}\left( t\right) +\beta _{2}\left[ \tilde{x}%
_{2}\left( t\right) \right] ^{2}}~,
\end{equation}%
\begin{equation}
\overset{\cdot }{\tilde{x}}_{2}\left( t\right) =\mathbf{i}\omega \tilde{x}%
_{2}\left( t\right) -\frac{\tilde{x}_{2}\left( t\right) +\alpha _{2}\tilde{x}%
_{1}\left( t\right) }{\beta _{1}\left[ \tilde{x}_{1}\left( t\right) \right]
^{2}+\left( \alpha _{1}\beta _{1}+\alpha _{2}\beta _{2}\right) \tilde{x}%
_{1}\left( t\right) \tilde{x}_{2}\left( t\right) +\beta _{2}\left[ \tilde{x}%
_{2}\left( t\right) \right] ^{2}}~,
\end{equation}%
\end{subequations}%
where $\mathbf{i=}\sqrt{-1}$ is the \textit{imaginary unit }and $\omega $ is
an \textit{arbitrary nonvanishing real} parameter, is \textit{isochronous:}
\textit{all} its nonsingular solutions are \textit{completely-periodic} in
the (\textit{real}) time variables $t,$ with a period which is $4$\textit{\
times} (or for a \textit{subset} of solutions only $2$\textit{\ times}) the
basic period $T=\pi /\left\vert \omega \right\vert $: see \cite{FC2008}. But
of course, due to the presence of the \textit{imaginary} unit $\mathbf{i}$
(see (\ref{isoiso})), this system (\ref{isoiso}) of ODEs lives in the
\textit{complex} world, namely its $2$ dependent variables $\tilde{x}%
_{1}\left( t\right) $ and $\tilde{x}_{2}\left( t\right) $ must, and its $4$
parameters $\alpha _{n}$ and $\beta _{n}$ ($n=1,2$) may, be considered to be
\textit{complex} numbers; so, in terms of \textit{real} variables and
\textit{real} parameters, it amounts to a system of $4$ \textit{nonlinear}
ODEs in $4$ \textit{real} dependent variables featuring $9$ \textit{%
arbitrary real }parameters. The fact that such a system---that the
interested reader might like to write out \textit{explicitly}---is \textit{%
isochronous} seems a \textit{remarkable} finding; possibly relevant in
\textit{applicative} contexts. $\blacksquare $

\textbf{Final remark}. The results reported above imply that the $2$
variables $x_{n}\left( t\right) $ or $\tilde{x}_{n}\left( t\right) $ evolve
in time as a \textit{linear} superposition with constant coefficients of $2$
simple functions of time: see (\ref{ODEs}). This is likely to motivate the
\textit{experts} on systems of nonlinear evolution equations to consider the
finding reported in this short communication to be rather \textit{trivial}.
On the other hand \textit{practitioners} might find that the systems of $2$
nonlinear systems of ODEs (\ref{ODEs}) or (\ref{isoiso}) describe
interesting phenomenologies. So it might be useful that the scientific
community become aware of this finding; which might possibly deserve to be
recorded in the website \textbf{EqWorld}. $\blacksquare $

 \strut\vfill

 \pagebreak

\subsection*{Acknowledgements}

The authors like to thank our colleague Robert
Conte for very fruitful discussions. FP would like to thank Payame Noor
University for financial support to this research.\label{lastpage}
\end{document}